\documentclass[a4paper,english]{smfart}

\usepackage{amsmath}
\usepackage{palatino}
\usepackage{amssymb}
\usepackage{mathrsfs}
\usepackage[all]{xy}
\usepackage{graphicx}
\usepackage{latexsym}
\usepackage{amsthm}
\usepackage{amscd}
\usepackage{mathrsfs}
\newtheorem{prop}{Proposition}
\newtheorem{theo}[prop]{Theorem}
\newtheorem{lemm}[prop]{Lemma}
\newtheorem{coro}[prop]{Corollary}
\theoremstyle{definition}

\newtheorem{rema}[prop]{Remark}

\numberwithin{prop}{section}

\def\q{{\bold Q}}
\def\z{{\bold Z}}
\def\qp{{\bold Q_p}}
\def\zp{{\bold Z_p}}
\def\qbar{\bar{\bold Q}}
\def\qpbar{\bar{\bold Q}_p}

\DeclareMathOperator\spf{Spf}

\begin{document}
\title{Lam\'e curves with bad reduction}
\author{Leonardo Zapponi}
\email{zapponi@math.jussieu.fr}
\urladdr{http://www.math.jussieu.fr/~zapponi}
\begin{abstract} Lam\'e curves are a particular class of elliptic curves (with a torsion point attached to them) which naturally arise when studying Lam\'e operators with finite monodromy. They can be realized as covers of the projective line unramified outside three points and can be defined over number fields. This paper investigates their $p$-adic properties. The main ingredient is formal/rigid geometry and in particular the use of $p$-adic theta functions. As a consequence, we can completely enumerate Lam\'e curves  with bad reduction (by giving, among the others, the $p$-adic valuation of their $j$-invariant) and describe the (local) Galois action.
\end{abstract}
\maketitle

\section*{Introduction} Following~\cite{Litcanu,Dahmen,Beukers}, the isomorphism classes of Lam\'e operators $L_n$ with finite monodromy bijectively correspond to a particular class of dessins d'enfants (isomorphism classes of covers of the projective line unramified outside three points). For operators $L_1$ with dihedral monodromy, these dessins are precisely those introduced in~\cite{Zapponi}. More precisely, from~\cite{Baldassarri}, the existence of such an operator is equivalent to the existence of a pair $(E,P)$, where $E$ is an elliptic curve and $P$ is a torsion point on it of order $n$ satisfying some extra conditions. In particular, there is a canonical degree $n$ cover $E\to\bold P^1$ which is unramified outside the points $\infty,0$ and $1$ and the operator corresponds to the associated dessin d'enfant. In this paper, the pair $(E,P)$ is called Lam\'e curve of order $n$. Remark that the curve $E$ can be defined over a number field and that that there exist finitely many Lam\'e curves of fixed order. This paper is mainly concerned with the study of their $p$-adic properties. The main ingredient is the formal/rigid geometry and especially the use of $p$-adic theta functions. This work is inspired by a previous investigation in~\cite{Zapponi}, where the problem was studied from a complex analytic point of view: it turned out that Lam\'e curves are closely related to (the zeroes of) a particular class of modular forms for the group $\Gamma_1(n)$. Here, we show that the same strategy can be adopted in the $p$-adic setting. This leads to a complete enumeration of Lam\'e curves with bad reduction. For example, we are now able to determine the $p$-adic valuation of their $j$-invariant and to describe the (local) Galois action.

The paper is organized as follows:

The first section is a brief introduction to Lam\'e curves and \S2 is a review of some known results in the archimdean case. The main result (Theorem 2.1) is an enumeration of Lam\'e curves of given order, arising from Grothendieck's theory of dessins d'enfants.

In \S3 we recall some well-known facts concerning Tate curves. More precisely, we desribe the $p$-adic uniformization of elliptic curves with bad reduction. We particularly focus on the construction of rational functions by using $p$-adic theta functions. (Proposition 3.1). We also include a 'cuspidal' parametrization of the modular curve $X_1(n)$ (Lemma 3.2), which is just an explicit construction of a moduli space (as rigid anaytic space) for pairs $(E,P)$ where $E$ is an elliptic curve with bad reduction and $P$ is a point of exact order $n$ on it.

In \S4 we introduce the functions $\varphi_\tau$, one for each cusp of $X_1(n)$. We then show (Proposition 4.1) that Lam\'e curves with bad reduction correspond to their zeroes. This is the $p$-adic analogue of the results in~\cite{Zapponi}.

In \S5 we obtain a complete enumeration of Lam\'e curves with bad reduction (Theorem 5.1) and we determine, among the others, the $p$-adic valuation of their $j$-invariants. This result leads to some rationality characterizations in \S6, where we determine the local field of moduli.

In the last section we give a detailed description of Lam\'e curves with bad reduction (at at least one prime, that is Lam\'e curves for which their $j$-invariant is not an algebraic integer) of order up to $11$. We also include a table for Lam\'e curves of order up to $20$.

\section{Lam\'e curves}
Let $E/\bold C$ be an elliptic curve and fix an element $P\in E[n]$ of exact order $n$ so that the divisor
$$D=n[P]-n[-P]$$
is principal, say $nD=(f)$, the rational function $f$ being defined up to a multiplicative constant. We say that the curve $E$ or, more correctly, that the pair $(E,P)$ is a {\it Lam\'e (elliptic) curve of order $n$} if the differential form $df$
has a (double) zero at the origin. We moreover assume that
$f(0_E)=1$, so that $f$ is uniquely determined. In this case the cover
$$E\to\bold P^1$$
induced by $f$ is unramified outside $\infty,0$ and $1$. A rigidity criterion of Weil asserts that the curve $E$ and the cover
itself can be defined over $\qbar$. In particular, the
$j$-invariant of $E$ is defined over a number field, this
condition being in fact equivalent to the arithmeticity of the cover.

\section{Review of the archimedean case}

Following~\cite{Dahmen,Litcanu,Zapponi}, the isomorphism classes of Lam\'e curves of given order over $\bold C$ can be enumerated via the Grothendieck's theory of dessins d'enfants (we refer to~\cite{Schneps} for a complete introduction to this subject). Skipping the intermediate steps, we directly give the final description.

In the following, we consider $3$-ples $(a,b,c)$, where $a,b$ and $c$ are positive integers. We moreover identify the  $3$-ples $(a,b,c)$ and $(b,c,a)$. The integer $n=a+b+c$ is the {\it degree} of $(a,b,c)$ and its {\it signature } is the canonical image of the integer $abc$ in $\bold F_2$, which is equal to $1$ if and only if $a,b$ and $c$ are odd. Finally, we say that $(a,b,c)$ is {\it primitive} if $\gcd(a,b,c)=1$

\begin{theo}\label{archimedian} Let $n$ be a positive integer. If $n$ is odd (resp. even) then there is a bijection between the set of $\bold C$-isomorphism classes of Lam\'e curves of order $n$ and the set of primitive $3$-ples of degree $n$ and signature $1$ (resp. of degree $n/2$ and signature $0$). 
\end{theo}

\begin{coro} There exist Lam\'e curves of any order $n\notin\{1,2,4,6\}$.
\end{coro}

\section{Tate curves and related topics}

\subsection{Tate curves} 
This is a short review of the theory of Tate
curves, we refer to~\cite{Silverman} for a complete introduction. Let $p$ be a prime number. For any $q\in\qpbar$ with $|q|<1$, the elliptic curve $E_q$ is defined as the quotient of $\qpbar^*$ with respect to the cyclic group $q^\bold Z$ generated by
$q$,
$$E_q=\qpbar^*/q^\bold Z$$
More precisely, $E_q$ is the rigid analytic space obtained from the annulus
$$\mathscr A_q=\{u\in\qpbar\,\,|\,\,|q|\leq|u|\leq 1\}$$
by identifying its boundaries via the multiplication by $q$. The curve $E_q$ is explicitely given by the affine equation
$$E_q:Y^2+XY=X^3+a_4X+a_6,$$
where $a_4=a_4(q)$ and $a_6=a_6(q)$ are power series which converge in the whole open unit disk. The discriminant of $E_q$ is simply given by the expression
$$\Delta(q)=q\prod_{n>0}(1-q^n)^{24}$$
and the $j$-invariant by the usual series
$$j(q)=q^{-1}+744+\dots\in\bold Z[\![q]\!].$$
The projection map $\qpbar^*\to E_q$, which is a homomorphism of groups, is then given in terms of power series ring;
its main advantage is that it commutes with the natural action of
$G_\qp=\mbox{Gal}(\qpbar/\qp)$ and, for any $p$-adic field
$K\subset\qpbar$ containing $q$, it induces an identification between
$K^*/q^\bold Z$ and $E_q(K)$.

Let $Y(1)=\mathcal H/\mbox{SL}_2(\bold Z)$ denote the coarse moduli space of elliptic curves over $\bold C$. It is a quasi-projective curve whose compactification $X(1)$ is isomorphic to the projective line. There is a unique cusp $\infty=X(1)\setminus Y(1)$, corresponding to the unique class of stable degenerated elliptic curve (that is, a projective line with two points identified). Both $Y(1)$ and $X(1)$ can be defined over $\bold Q$ and, after a base extension to $\qp$, they inherit a natural structure of rigid analytic space. More precisely, there is a canonical smooth affine (formal) model of $Y(1)$ over $\zp$, namely
$$\mathcal Y(1)=\spf(\zp\{j\})$$
and, as a rigid analytic space, $Y(1)$ is its generic fiber. The complement
$$\Bbb D_\infty=X(1)-Y(1)=\mbox{Spm}(K\langle j^{-1}\rangle)$$
is an open unit rigid disk. For any $p$-adic field $K$,  the elements of $\Bbb D_\infty(K)$ correspond to the $K$-rational points of $X(1)$ specializing to $\infty$, that is to the $\qpbar$-isomorphism classes of ellitpic curves $E/K$ with bad reduction. Remark that since the above power series expansion of $j$ induces an automorphism of the unit open rigid disk, we obtain
$$\Bbb D_\infty\cong\mbox{Spm}(K\langle q\rangle).$$
In other words, $q$ is a global parameter for elliptic curves with bad reduction over $K$, i.e. a $\qpbar$-isomorphism class of elliptic curves $E/K$ corresponds to a unique $q\in K^*$ with $|q|<1$. We could have choosen $j$ as parameter, which is somehow more natural, but as we will see in the rest of the paper, $q$ is more convenient.
 
\subsection{Theta functions and rational functions}
As in the archimedian case, we can express rational funtions on $E_q$
by using the normalized theta function
$$\theta(u)=(1-u)\prod_{n>0}{(1-uq^n)(1-u^{-1}q^n)\over(1-q^n)^2}$$
which converges for any $u\in\qpbar^*$ and any $q\in\qpbar$ with
$|q|<1$. Indeed, let
$$D=\sum_ie_iP_i$$
be a divisor on $E_q$ and fix $u_1,\dots,u_n\in\qpbar^*$ lying above
$P_1,\dots,P_n$ with respect to the canonical projection $\qpbar^*\to
E_q$, which is a homomorphism of abelian groups. Then $D$ is principal if and only if there exists an integer
$\nu$ such that
$$\prod_{i=1}^nu_i^{e_i}=q^\nu.$$
The following is a classical result which holds for any complete field.

\begin{prop} The notation and hypothesis being as above, the (pull-back in $\qpbar^*$ of a) rational function $f$ on $E$ having $D$ as divisor is given by
$$f(u)=u^{-\nu}\prod_{i=1}^n\theta(u_i^{-1}u)^{e_i}.$$
\end{prop}

\subsection{Torsion points}

Let $E_q$ be the tate curve corresponding to $q\in\qpbar$ with
$|q|<1$. Any $n$-torsion point $P$ of
$E_q$ has a unique representative $u\in\qpbar^*$ with $|q|<|u|\leq1$,
which can be written as
$$u=\xi\rho^b,$$
where $\xi$ is a $n$-th root of unity, $\rho$ is a $n$-th root of $q$ and $0\leq b<n$ is an integer. Setting $d=\gcd(n,b)$, $n'=n/d$ and $b'=b/d$, we find
$$u^{n'}=\zeta q^{b'},$$
where $\zeta$ is a primitive $d$-th root of unity (as soon as $P$ has exact order $n$). In the following, we say that $P$ is {\it of type} $\tau=(b,\zeta)$. We moreover set
$$K_\tau=\qp(\zeta)=\qp(\boldsymbol\mu_d)$$
and denote by $R_\tau$ its ring of integers. There is a natural action of $G_\qp$ on the set of types, namely
$$\sigma(b,\zeta)=(b,\sigma(\zeta)).$$
In this case, if a $n$-torsion point of $E_q$ is of type $\tau$ then the point $\sigma(P)$ of $^\sigma E_q=E_{\sigma(q)}$ is of type $\sigma(\tau)$ and $K_\tau$ is the 'field of definition' of $\tau$.

\subsection{Cuspidal parametrization of $X_1(n)$}
Let $Y_1(n)=\mathcal H/\Gamma_1(n)$ be the coarse moduli space of pointed elliptic curves $(E,P)$ over $\bold C$, where $P$ is a point of exact order
$n$. As in the case of $Y(1)$, it is a smooth quasi-projective curve; we denote by $X_1(n)$ its compactification. Let $S=X_1(n)-Y_1(n)$ be the set of cusps. Both of these curves can be defined over $\bold Q$ and, over  $\qp$, they have a natural structure of rigid analytic space. There is a natural forgetful map
$$X_1(n)\to X(1)$$
which is in fact a finite cover, unramified outside
$\infty,0$ and $1728$. The pull-back
$$\Bbb D_n=\Bbb D_\infty\times_{X(1)}X_1(n)=\bigcup_{P\in S}\Bbb D_P$$
is the disjoint union of open unit rigi disks, one for each cusp of $X_1(n)$. We now describe how they can be explicitely constructed. As we have seen in the previous paragraph, each $n$-torsion point $P$ on $E_q$ is of a certain
type $\tau=(b,\zeta)$. Since we are working up to isomorphism, and since
$(E_q,P)$ is isomorphic to $(E_q,-P)$, we can assume that $b\leq n/2$
(remark that $(E_q,-P)$ is of type $(n-b,\zeta^{-1})$). Moreover, for $b=0$, we have to identitfy the curves of type $(0,\zeta)$ with the curves of type $(0,\zeta^{-1})$. With these conventions, a type uniquely corresponds to a cusp of $X_1(n)$. For any type $\tau$, consider the rigid open unit disk
$$\Bbb D_\tau=\mbox{Spm}(\qp\langle\rho\rangle).$$

\begin{lemm}\label{cuspidal} There is a natural bijection between $\Bbb D_\tau(\qpbar)-0$ and the set of $\qpbar$-isomorphism classes of pointed elliptic curves $(E,P)$, where $P\in E[n]$ is of type $\tau$. Moreover, this correspondence is compatible with the action of $G_\qp$.
\end{lemm}

\begin{proof} For any $b<n$, following the notation of the previous paragraph, fix once for all two positive integers $b'',n''$ such that
$$n'n''-b'b''=1.$$
Let $\tau=(b,\zeta)$ be a type. Given $\rho\in\Bbb D_\tau(\qpbar)-0$, set
$$\left\{\aligned
&q=\zeta^{b''}\rho^{n'},\\
&u=\zeta^{n''}\rho^{b'}.
\endaligned\right.$$
The pair $(E,P)$ corresponding to $\rho$ is defined as follows: $E=E_q$ is the Tate curve associated to $q$ and $P$ is the $n$-torsion point corresponding to $u\in\qpbar^*$, which is of type $\tau$.

Conversely, given a Tate curve $E_q$ and a $n$-torsion point $P$ on it
of type $\tau$, corresponding to $u\in\qpbar^*$ with $|q|<|u|\leq1$ then $\rho$ is given by the relation
$$\rho=u^{-b''}q^{n''}.$$
\end{proof}

\begin{rema} Following the notation of the previous paragraph, the  field $K_\tau$ is precisely the field of definition of the cusp of $X_1(n)$ corresponding to the type $\tau$.\end{rema}

\section{The functions $\varphi_\tau$}
We now start the study Lam\'e curves over $\qpbar$ of fixed order $n$
with bad reduction. As in \S3.3, for any $b<n$ set $d=\gcd(n,b)$, $n'=n/d$, $b'=b/d$ and fix two positive integers $n'',b''$ such that
$$n'n''-b'b''=1.$$
Given a type $\tau=(b,\zeta)$, the power series $\varphi_\tau$ defined by
$$\varphi_\tau=\frac{n-2b}{2n}+\sum_{m\geq0}\frac{\zeta^{mb''+n''}\rho^{mn'+b'}} {1-\zeta^{mb''+n''}\rho^{mn'+b'}}-\sum_{m>0}\frac{\zeta^{mb''-n''}\rho^{mn'-b'}}{1-\zeta^{mb''-n''}\rho^{mn'-b'}}\in K_\tau[\![\rho]\!]$$
converges for any $\rho\in\Bbb D_\tau(\qpbar)$.

\begin{prop}\label{prop1} Given $\rho\in\Bbb D_\tau(\qpbar)-0$, the corresponding pair $(E,P)$ is a Lam\'e curve of order $n$ if and only if $\varphi_\tau(\rho)=0$.
\end{prop}

\begin{proof}
Denote by $E=E_q$ be the Tate curve associated to $q=\zeta^{b''}\rho^{n'}$, the torsion point $P$ corresponding to the element $v=\zeta^{n''}\rho^{b'}\in\qpbar^*$ (cf. the previous section). We know from \S3.2 that a rational function $f\in\qpbar(E)$ such that $(f)=n[P]-n[0_E]$ can uniquely be written as
$$f(u)=cu^{-2b}{\theta(v^{-1}u)^n\over\theta(vu)^n},$$
with $c\in\qpbar^*$. Taking logarithmic derivatives, we get the identity
$$\frac{df}f=\psi_\tau(u)\frac{du}u,$$
where the function $\psi_\tau$ is defined by
$$\aligned
\psi_\tau(u)=&-2b-n\sum_{m\geq0}\frac{v^{-1}uq^m}{1-v^{-1}uq^m}+ n\sum_{m>0}\frac{vu^{-1}q^m}{1-vu^{-1}q^m}+\\
&+n\sum_{m\geq0}\frac{vuq^m}{1-vuq^m}- n\sum_{m>0}\frac{v^{-1}u^{-1}q^m}{1-v^{-1}u^{-1}q^m}
\endaligned$$
and converges for any $u\in\qpbar^*$ not belonging to $vq^\z\cup v^{-1}q^\z$. Now, the pair $(E,P)$ is a Lam\'e curve of order $n$ if and only if $df$ has a zero at $0_E$, which can be restated as $\psi_\tau(1)=0$. The above expression of $\psi_\tau$ leads to the identity
$$\aligned
\psi_\tau(1)&=-2b-\frac{nv^{-1}}{1-v^{-1}}+\frac{nv}{1-v}-2n\sum_{m>0}\frac{v^{-1}q^m}{1-v^{-1}q^m} +2n\sum_{m>0}\frac{vq^m}{1-vq^m}=\\
&=n-2b+2n\sum_{m\geq0}\frac{vq^m}{1-vq^m}-2n\sum_{m>0}\frac{v^{-1}q^m}{1-v^{-1}q^m}.
\endaligned$$
If we express $q$ and $u$ in terms of $\rho$ and $\zeta$, we finally obtain
$$\aligned
\psi_\tau(1)&=n-2b+2n\sum_{m\geq0}\frac{\zeta^{mb''+n''}\rho^{mn'+b'}}{1-\zeta^{mb''+n''}\rho^{mn'+b'}}-2n\sum_{m>0}\frac{\zeta^{mb''-n''}\rho^{mn'-b'}}{1-\zeta^{mb''-n''}\rho^{mn'-b'}}=\\
&=2n\varphi_\tau(\rho).
\endaligned$$
The condition $\psi_\tau(1)=0$ is now clearly equivalent to $\varphi_\tau(\rho)=0$.
\end{proof}

\section{Existence}
We now give an existence criterion for Lam\'e curves with bad reduction of given type.

\begin{theo}\label{th2} Given an integer $n>2$ and a corresponding type $\tau=(b,\zeta)$ with $b\leq n/2$, the following conditions are equivalent:
\begin{enumerate}
\item There exists a Lam\'e curve of order $n$ over $\qpbar$ with bad reduction of type $\tau$.
\item The inequality $|n-2b|<|2n|$ holds.
\end{enumerate}
If one of these conditions is fulfilled then there exists exactly $b'=b/\gcd(n,b)$ isomorphism classes of Lam\'e curves of order $n$ with bad reduction of type $\tau$, which have a natural structure of $\mu_{b'}$-torsor. Moreover, if $j$ denotes the $j$-invariant of such  a curve then we have the identity
$$|j|=\left|\frac{2n}{n-2b}\right|^{n/b}.$$
\end{theo}

\begin{proof} Let $\tau=(b,\zeta)$ be a type. From now on, we assume that $b\leq n/2$. Since we are working up to isomorphism, this is not a restriction, cf.~\S3.4. We know from Proposition~\ref{prop1} that the isomorphism classes of Lam\'e curve of type $\tau$ correspond to the zeroes of $\varphi_\tau$ in the  pointed rigid unit disk $\Bbb D_\tau-0$. If $b=0$, that is if $\tau=(0,\zeta)$ with $\zeta$ a primitive $n$-th root of unity, then we have $d=n$, $n'=1$, $b'=0$; we then take $n''=1$ and $b''=0$ and we find a power series expansion of the form
$$\varphi_\tau=(1+\zeta)\left(\frac1{2(1-\zeta)}+\rho f_\tau\right),$$
with $f_\tau\in R_\tau[\![\rho]\!]$, so that $\varphi_\tau$ cannot have zeroes in the open unit rigid disk. Remark that $1+\zeta$ vanishes if and only if $n=2$. This case can be excluded, since following Theorem 2.1, there are no Lam\'e curves of order $2$. For $b=n/2$, we find $\tau=(n/2,\zeta)$, where $\zeta$ is a primitive $b$-th root of unity. We have $d=n/2$, $n'=2$ and $b'=1$; taking $n''=1$ and $b''=1$ we then find a power series expansion of the form
$$\varphi_\tau=(1-\zeta)\rho f_\tau,$$
with $f_\tau\in R_\tau[\![\rho]\!]^\times$. In particular, $\rho=0$ is its unique zero in the open unit rigid disk. Once again, $1-\zeta$ vanishes if and only if $n=2$, which is excluded. We can therefore assume that $b\neq0,n/2$. In this case, the function $\varphi_\tau$ has a power series expansion of the form
$$\varphi_\tau=c_0+\rho^{b'}f_\tau,$$
where $c_0=\frac{n-2b}{2n}$ and $f_\tau$ is a unit in $R_\tau[\![\rho]\!]$. In particular, if $\varphi_\tau(\rho)=0$ then we obtain
$$|c_0|=|\rho^{b'}f_\tau(\rho)|=|\rho|^{b'}<1.$$
This proves that the inequality in condition (2) is necessary. Moreover, since $q=\zeta^{b''}\rho^{n'}$ and since $|j|=|q|^{-1}$, the last equality follows.

Suppose now that the inequality in condition (2) is fulfilled, with $b\neq0,n/2$. Remark that in this case $p$ does not divides $b'$. Setting
$$\tilde\varphi_\tau=\varphi_\tau-\frac{n-2b}{2n},$$
we then find
$$\tilde\varphi_\tau=\zeta^{n''}\rho^{b'}+u_1\rho^{b'+1}+\cdots=\rho^{b'}f_\tau\in R_\tau[\![\rho]\!],$$
where $f_\tau$ is a unit in $R_\tau[\![\rho]\!]$. In particular, there exists a power series
$$g_\tau=\rho+v_2\rho^2+\cdots\in R_\tau[\![\rho]\!]$$
inducing an isomorphism $g_\tau:\Bbb D_\tau\to\Bbb D_\tau$ defined over $K_\tau$ such that
$$\tilde\varphi_\tau\circ g_\tau(\rho)=\zeta^{n''}\rho^{b'}.$$
In this case, the zeroes of $\varphi_\tau$ are the images under $g_\tau$ of the $b'$-th roots of
$$\frac{2b-n}{2\zeta^{n''}n}$$
and the $\mu_{b'}$-torsor structure naturally arises from this last construction.
\end{proof}

\begin{coro}\label{cor1} Let $p$ be a prime number and consider an integer $n=p^rm$ with $p\!\not |m$. If $m<p$ then any Lam\'e curve of order $n$ over $\qpbar$ has good reduction.
\end{coro}

\begin{proof} We know from Theorem~\ref{th2} that the existence of a Lam\'e curve of order $n$ over $\qpbar$ with bad reduction is equivalent to the existence of a positive integer $b<n/2$ such that
$$v(n-2b)>v(2n)\geq r.$$
Here $v$ denotes the $p$-adic valuation, normalized by the relation $v(p)=1$. In particular, the above inequality is equivalent to the existence of an integer $n'<n$ divisible by $p^{r+1}$, which is impossible if $m<p$.
\end{proof}

This last result is not optimal, since it does not completely describes the orders for which any Lam\'e curve over $\qpbar$ has good reduction. For example, we have the following

\begin{coro}\label{cor2} Any Lam\'e curve of odd order $n$ over $\qbar_2$ has good reduction.
\end{coro}

\begin{proof} Indeed, for any positive integer $b<n/2$, the integer $n-2b$ is odd and thus its $2$-adic valuation cannot be greater than the valuation of $2n$.
\end{proof}

\section{Rationality questions, local Galois action}

We now investigate the behaviour of the local field of moduli of a Lam\'e curve $(E,P)$ over $\qpbar$, that is, the subfield of $\qpbar$ fixed by the elements $\sigma\in G_\qp$ such that $(^\sigma E,^\sigma P)$ is isomorphic to $(E,P)$. For any type $\tau=(b,\zeta)$  with $2b\neq n$, set $d=\gcd(n,b), n'=n/d$ and $b'=b/d$ and consider the element
$$\alpha_\tau={2b-n\over2\zeta^mn}\in K_\tau^\times/(K_\tau^\times)^{b'}.$$
where the integer $m$ is an inverse of $n'$ modulo $b'$.

\begin{prop}\label{moduli} Let $K$ denote the local field of moduli of a Lam\'e curve of order $n$ and type $\tau=(b,\zeta)$ over $\qpbar$, with $b\leq n/2$. We then have
$$K=K_\tau(\theta),$$
where $\theta$ is a $b'$-th root of $\alpha_\tau$.
\end{prop}

\begin{proof} Remark first of all that the field of moduli (which is also a field of definition) of an element $\rho\in\Bbb D_\tau(\qpbar)$ is simply $K_\tau(\rho)$. In particular, it contains $K_\tau$. We follow the notation of the proof of Theorem~\ref{th2}, where we proved that an element $\rho\in\Bbb D_\tau(\qpbar)$ corresponds to a Lam\'e curve if and only if it is the image under the automorphism $g_\tau$ of a $b'$-th root $\alpha$ of
$${2b-n\over2\zeta^{n''}n}.$$
We then obtain $K=K_\tau(g_\tau(\alpha))$. We want to prove that $K=K_\tau(\alpha)$. But this directly follows from the fact that $g_\tau$ is defined over $K_\tau$. Indeed, for any $\sigma\in\mbox{Gal}(\qpbar/K_\tau)$, we obtain
$$\sigma(\alpha)=\alpha\Leftrightarrow g_\tau(\sigma(\alpha))=g_\tau(\alpha)\Leftrightarrow\sigma(g_\tau(\alpha))=g_\tau(\alpha).$$
Finally, the identity
$$n'n''-b'b''=1$$
implies that $n''$ is an inverse of $n'$ modulo $b'$. In particular $b'$ divides $m-n''$, so that $(2b-n)/2\zeta^mn$ and $(2b-n)/2\zeta^{n''}n$ define the same element in $K_\tau^\times/(K_\tau^\times)^{b'}$ and thus there exists a $\theta$ such that $K_\tau(\alpha)=K_\tau(\theta)$.

\end{proof}

\begin{coro} The notation being as in Proposition~\ref{moduli}, the extension $K/K_\tau$ is tamely ramified.
\end{coro}

\begin{proof} We know from Theorem~\ref{th2} that
$$v(n-2b)>v(2n)\geq v(n).$$
In particular, if $v(2b)>v(n)$ we obtain
$$v(n)=v(n-2b)>v(n),$$
which is impossible. We then have
$$v(b)\leq v(2b)\leq v(n)$$
and thus $p$ does not divides $b'$ so that the extension $K/K_\tau$ is tamely ramified.
\end{proof}

\begin{coro} Suppose that $p$ is odd (resp. that $p=2$) and that it divides $n-2$ (resp. that $n$ is congruent to $2$ modulo $8$). Then there exists a Lam\'e curve with bad reduction defined over $\qp$.
\end{coro}

\begin{proof} Consider the type $\tau=(1,1)$. Since $v(n-2)>v(2n)$, Theorem~\ref{th2} asserts that there exist exactly one Lam\'e curve over  $\qpbar$ with bad reduction of type $\tau$. In this case, Proposition~\ref{moduli} asserts that its field of definition is $K_\tau=\qp$.
\end{proof}

\section{Examples}

We close this paper with the complete enumeration of Lam\'e curves of order up to $11$ having bad reduction at some prime, i.e. for which the $j$-invariant is not an algebraic integer.

\subsection{Order $5$}

Up to isomorphism, there exists a unique Lam\'e curve of order $5$ over $\qbar$, corresponding to the $3$-ple $(1,1,3)$ (cf.~\S2). In particular, it is automatically defined over $\q$ and its $j$-invariant is a rational number. In~\cite{Zapponi}, we found
$$j={20480\over243}={2^{12}\cdot5\over3^5},$$
so that $p=3$ is the only prime of bad reduction. We can now directly check the results of the previous sections. First of all, Corollary~\ref{cor1} and Corollary~\ref{cor1} assert that $p=3$ is the only possible prime of bad reduction. The only possible value of $b\leq n/2$ is $b=1$. In particular, the curve is of type $\tau=(1,1)$ and the assosiated power series $\varphi_\tau$ is given by
$$\varphi_\tau=\frac3{10}+\sum_{m\geq0}{\rho^{5m+1}\over1-\rho^{5m+1}}-\sum_{m>0}{\rho^{5m-1}\over1-\rho^{5m-1}}.$$
Finally, Theorem~\ref{th2} asserts that there exists a unique isomorphism class of Lam\'e  curves with reduction of type $\tau$ at $p=3$ and that the $3$-adic valuation of $j$ is equal to $-n/b=-5$.

\subsection{Order $7$}

There exist two isomorphism classes of Lam\'e curve of order $7$ over $\qbar$, corresponding to the $3$-ples $(1,1,5)$ and $(1,3,3)$. They are both conjugated under $G_\q$ and thus they are defined over a quadratic extension $K$ of $\q$. A direct computation gives $K=\q(\sqrt{21})$. The only poosible primes of bad reduction are $3$ and $5$.

For $p=3$ we find $b=2$ and following Theorem~\ref{th2} there exist two isomorphism classes of Lam\'e curves of type $\tau=(2,1)$, associated to the power series
$$\varphi_\tau=\frac3{14}+\sum_{m\geq0}{\rho^{7m+2}\over1-\rho^{7m+2}}-\sum_{m>0}{\rho^{7m-2}\over1-\rho^{7m-2}}.$$
This implies that any Lam\'e curve of order $7$ have bad reduction at (the unique prime of $K$ above) $3$. Moreover, the valuation of its $j$-invariant is equal to $-n/b=-7/2$.

For $p=5$ we find $b=1$ so that there exists a unique isomorphism class of Lam\'e curves of type $\tau=(1,1)$, associated to the power series
$$\varphi_\tau=\frac5{14}+\sum_{m\geq0}{\rho^{7m+1}\over1-\rho^{7m+1}}-\sum_{m>0}{\rho^{7m-1}\over1-\rho^{7m-1}}.$$
The valuation of the associated $j$-invariant is equal to $-n/b=-7$. This results have the following arithmetic consequences: the prime $p=5$ splits in $K$, say $5\mathcal O_K=\frak p_1\frak p_2$, and there exists a unique isomorphism class for which the $\frak p_1$-valuation (resp. the $\frak p_2$-valuation) of the associated $j$-invariant is equal to $-7$ (resp. to $7$).

\subsection{Order $8$}

There exists a unique isomorphism class of Lam\'e curves of order $8$ over $\qbar$, corresponding to the $3$-ple $(1,1,2)$. Following~\cite{Zapponi}, the corresponding $j$-invariant is given by
$$j={210646\over6561}={210646\over3^8}$$
The situation is similar to the order $5$ case: the curve is of type $\tau=(1,1)$, the assosiated power series $\varphi_\tau$ is given by
$$\varphi_\tau=\frac38+\sum_{m\geq0}{\rho^{8m+1}\over1-\rho^{8m+1}}-\sum_{m>0}{\rho^{8m-1}\over1-\rho^{8m-1}}$$
and the $3$-adic valuation of $j$ is equal to $-n/b=-8$.

\subsection{Order $9$}

There are three isomorphism classes of Lam\'e curves of order $9$ over $\qbar$, corresponding to the $3$-ples $(1,1,7),(1,3,5)$ and $(1,5,3)$. The only possible primes of bad reduction are $5$ and $7$.

For $p=5$ we find $b=2$ and thus there exist two isomorphism classes of Lam\'e curves of type $\tau=(2,1)$, associated to the power series
$$\varphi_\tau=\frac5{18}+\sum_{m\geq0}{\rho^{9m+2}\over1-\rho^{9m+2}}-\sum_{m>0}{\rho^{9m-2}\over1-\rho^{9m-2}}$$
and the valuation of their $j$-invariants is equal to $-n/b=-9/2$.

For $p=7$ we find $b=1$ so that there exists a unique isomorphism class of Lam\'e curves of type $\tau=(1,1)$, associated to the power series
$$\varphi_\tau=\frac7{18}+\sum_{m\geq0}{\rho^{9m+1}\over1-\rho^{9m+1}}-\sum_{m>0}{\rho^{9m-1}\over1-\rho^{9m-1}}.$$
The valuation of the associated $j$-invariant is equal to $-n/b=-9$.

\subsection{Order $10$}

There is a unique isomorphism class of Lam\'e curves of order $10$ over $\qbar$, corresponding to the $3$-ple $(1,2,2)$. The associated $j$-invariant being a rational number, we now determine its denominator. The only primes of bad reduction are $2$ and $3$.

For $p=2$ we find $b=1$ so that the curve is of type $\tau=(1,1)$, with assotiated power series
$$\varphi_\tau=\frac25+\sum_{m\geq0}{\rho^{10m+1}\over1-\rho^{10m+1}}-\sum_{m>0}{\rho^{10m-1}\over1-\rho^{10m-1}}.$$
The valuation of the $j$-invariant is then equal to $-n/b=-10$.

For $p=3$ we find $b=2, d=2, b'=1$ and $n'=5$, so that the curve is of type $\tau=(2,-1)$. By setting $n''=1$ and $b''=0$ we obtain the power series
$$\varphi_\tau=\frac3{10}-\sum_{m\geq0}{\rho^{5m+1}\over1+\rho^{5m+1}}+\sum_{m>0}{\rho^{5m-1}\over1+\rho^{5m-1}}$$
and the valuation of the $j$-invariant is equal to $-n/b=-5$.

We finally deduce that the denominator of the $j$-invariant is equal to $248832=2^{10}\cdot3^5$. Remark that this is the first case where the bad reduction at $2$ occurs.

\subsection{Order $11$}

There are five isomorphism classes of Lam\'e curves of order $11$ over $\qbar$, corresponding to the $3$-ples $(1,1,9),(1,3,7),(1,7,3),(1,5,5)$ and $(3,3,5)$. The only possible primes of bad reduction are $3,5$ and $7$.

For $p=3$ we find $b=1$ or $b=4$. In the first case, there exists a unique isomorphism class with bad reduction of type $\tau=(1,1)$ and assotiated power series
$$\varphi_\tau=\frac9{22}+\sum_{m\geq0}{\rho^{11m+1}\over1-\rho^{11m+1}}-\sum_{m>0}{\rho^{11m-1}\over1-\rho^{11m-1}}.$$
The valuation of the $j$-invariant is equal to $-2n/b=-22$. For $b=4$ we obtain four isomorphism classes with bad reduction of type $(4,1)$ and assotiated power series
$$\varphi_\tau=\frac3{22}+\sum_{m\geq0}{\rho^{11m+4}\over1-\rho^{11m+4}}-\sum_{m>0}{\rho^{11m-4}\over1-\rho^{11m-4}}.$$
The valuation of their $j$-invariants is equal to $-n/b=-11/4$. In particular, we deduce that any Lam\'e curve of order $11$ over $\bar{\bold Q}_3$ has bad reduction.

For $p=5$ we find $b=3$ so that there exist three isomorphism classes with bad reduction, of type $\tau=(3,1)$ and assotiated power series
$$\varphi_\tau=\frac5{22}+\sum_{m\geq0}{\rho^{11m+3}\over1-\rho^{11m+3}}-\sum_{m>0}{\rho^{11m-3}\over1-\rho^{11m-3}}.$$
The valuation of their $j$-invariants is equal to $-n/b=-11/3$.

For $p=7$ we find $b=2$ so that there exist two isomorphism classes with bad reduction, of type $\tau=(2,1)$ and assotiated power series
$$\varphi_\tau=\frac7{22}+\sum_{m\geq0}{\rho^{11m+2}\over1-\rho^{11m+2}}-\sum_{m>0}{\rho^{11m-2}\over1-\rho^{11m-2}}.$$
The valuation of their $j$-invariants is equal to $-n/b=-11/2$.
\begin{center}
\begin{tabular}{|c|c|c|c|c|c|c|c|}
\hline
$n$ & \# & $p$ & \# & $\tau$ & $v(j)$ & $[K:\qp]$\\
\hline
$12$ & $3$ & $5$ & $1$ & $(1,1)$ & $-12$ & $1$\\
\hline
$13$ & $7$ & $3$ & $2$ & $(2,1)$ & $-13$ & $2$\\
& & & $5$ & $(5,1)$ & $-13/5$ & $5$\\
& & $5$ & $4$ & $(4,1)$ & $-13/4$ & $4$\\
& & $7$ & $3$ & $(3,1)$ & $-13/3$ & $3$\\
& & $11$ & $1$ & $(1,1)$ & $-13$ & $1$\\
\hline
$14$ & $3$ & $2$ & $3$ & $(3,1)$ & $-14/3$ & $3$\\
& & $3$ & $1$ & $(1,1)$ & $-14$ & $1$\\
& & & $2$ & $(4,-1)$ & $-7/2$ & $2$\\
& & $5$ & $1$ & $(2,-1)$ & $-7$ & $1$\\
\hline
$15$ & $8$ & $3$ & $2$ & $(3,\zeta_3)$ & $-5$ & $2$\\
& & $7$ & $4$ & $(4,1)$ & $-15/4$ & $4$\\
& & $11$ & $2$ & $(2,1)$ & $-15/2$ & $2$\\
& & $13$ & $1$ & $(1,1)$ & $-15$ & $1$\\
\hline
$16$ & $6$ & $3$ & $1$ & $(2,-1)$ & $-8$ & $1$\\
& & & $5$ & $(5,1)$ & $-16/5$ & $5$\\
& & $5$ & $3$ & $(3,1)$ & $-16/3$ & $3$\\
& & $7$ & $1$ & $(1,1)$ & $-16$ & $1$\\
\hline
$17$ & $12$ & $3$ & $1$ & $(1,1)$ & $-17$ & $1$\\
& & & $4$ & $(4,1)$ & $-17/2$ & $4$\\
& & & $7$ & $(7,1)$ & $-17/7$ & $7$\\
& & $5$ & $1$ & $(1,1)$ & $-17$ & $1$\\
& & & $6$ & $(6,1)$ & $-17/6$ & $6$\\
& & $7$ & $5$ & $(5,1)$ & $-17/5$ & $5$\\
& & $11$ & $3$ & $(3,1)$ & $-17/3$ & $3$\\
& & $13$ & $2$ & $(2,1)$ & $-17/2$ & $2$\\
\hline
\end{tabular}
\end{center}
\begin{center}
\begin{tabular}{|c|c|c|c|c|c|c|c|}
\hline
$n$ & \# & $p$ & \# & $\tau$ & $v(j)$ & $[K:\qp]$\\
\hline
$18$ & $6$ & $2$ & $1$ & $(1,1)$ & $-36$ & $1$\\
& & & $5$ & $(5,1)$ & $-18/5$ & $5$\\
& & $5$ & $2$ & $(4,-1)$ & $-9/2$ & $2$\\
& & $7$ & $1$ & $(2,-1)$ & $-9$ & $1$\\
\hline
$19$ & $15$ & $3$ & $2$ & $(2,1)$ & $-19/2$ & $2$\\
& & & $5$ & $(5,1)$ & $-38/5$ & $5$\\
& & & $8$ & $(8,1)$ & $-19/8$ & $8$\\
& & $5$ & $2$ & $(2,1)$ & $-19/2$ & $2$\\
& & & $7$ & $(7,1)$ & $-19/7$ & $7$\\
& & $7$ & $6$ & $(6,1)$ & $-19/6$ & $6$\\
& & $11$ & $4$ & $(4,1)$ & $-19/4$ & $4$\\
& & $13$ & $3$ & $(3,1)$ & $-19/3$ & $3$\\
& & $17$ & $1$ & $(1,1)$ & $-19$ & $1$\\
\hline
$20$ & $10$ & $2$ & $1$ & $(2,-1)$ & $-10$ & $1$\\
& & $3$ & $1$ & $(1,1)$ & $-40$ & $1$\\
& & & $2$ & $(4,\zeta_4)$ & $-5$ & $2$\\
& & & $7$ & $(7,1)$ & $-20/7$ & $7$\\
& & $7$ & $3$ & $(3,1)$ & $-20/3$ & $3$\\
\hline
\end{tabular}
\end{center}

\end{document}